\documentclass[amstex,12pt, amssymb]{article}

\usepackage{mathtext}
\usepackage[cp1251]{inputenc}
\usepackage[T2A]{fontenc}
\usepackage[dvips]{graphicx}
\usepackage{amsmath}
\usepackage{amssymb}
\usepackage{amsxtra}
\usepackage{latexsym}
\usepackage{ifthen}

\newcounter{lemma}[section]

\newcounter{corollary}[section]

\newcounter{remark}[section]

\newcounter{theorem}[section]

\newcounter{proposition}[section]

\newcounter{example}

\numberwithin{equation}{section}

\begin{document}

\markboth{V. DESYATKA, A.~HALYTS'KA, E.~SEVOST'YANOV}{\centerline{ON
BOUNDARY EXTENSION ...}}

\def\cc{\setcounter{equation}{0}
\setcounter{figure}{0}\setcounter{table}{0}}

\overfullrule=0pt

%\normalsize\large

\author{VICTORIA DESYATKA, ALINA HALYTS'KA, EVGENY SEVOST'YANOV}

\title{
{\bf ON BOUNDARY EXTENSION OF UNCLOSED ORLYCZ-SOBOLEV MAPPINGS}}

\date{\today}
\maketitle

%\large
\begin{abstract}
This paper is devoted to the study of the boundary behavior of
Orlicz-Sobolev classes that may not preserve the boundary under
mapping. Under certain conditions, we show that these mappings have
a continuous boundary extension of the definition domain.
\end{abstract}

\bigskip
{\bf 2010 Mathematics Subject Classification: Primary 30C65;
Secondary 31A15, 31B25}

\section{Introduction}

The present paper is devoted to the study of the boundary behavior
of mappings with bounded and finite distortion. Various authors have
obtained many results on this topic, see, e.g., \cite{Cr},
\cite{MRSY$_1$}--\cite{MRSY$_2$}, \cite{Na$_1$}--\cite{Na$_2$} and
\cite{Vu}. In particular, relatively recently some authors have
studied the Orlicz-Sobolev classes and their connection with the
distortion of moduli of families of paths and surfaces, see e.g.
\cite{RS$_2$}, \cite{KR}, \cite{RSS} and \cite{SSP}. We also
obtained some results on the boundary behavior of these mappings,
see~\cite{Sev$_1$} and \cite{Sev$_2$}. In~\cite{Sev$_2$}, we have
established a continuous boundary extension of open discrete
Orlicz-Sobolev mappings which are closed (equivalently, boundary
preserving). Similar studies have not been carried out for
non-closed mappings, and our goal now is to investigate their. Some
steps in this direction have already been made in~\cite{DS}. Note
that, the connection between Orlicz-Sobolev classes and estimates of
the distortion of the modulus of families of paths has been
established only in some special cases, say, for homeomorphisms, or,
more generally, open, discrete and closed mappings (see e.g.
\cite[Lemma~2, Theorem~5]{Sev$_2$}). The results which we obtain
here are new in relation to the results mentioned above.

\medskip
Recall some definitions. Let $U$ be an open set in ${\Bbb R}^n.$ In
what follows, $C^k_0(U)$ denotes the space of functions $u:U
\rightarrow {\Bbb R} $ with a compact support in $U,$ having $k$
partial derivatives with respect to any variable that are continuous
in $U.$ Let $U$ be an open set, $U\subset{\Bbb R}^n,$ and let
$u:U\rightarrow {\Bbb R}$ be some function, $u \in L_{\rm
loc}^{\,1}(U).$ Suppose that, there is a function $v\in L_{\rm
loc}^{\,1}(U)$ such that
$$\int\limits_U \frac{\partial \varphi}{\partial x_i}(x)u(x)\,dm(x)=
-\int\limits_U \varphi(x)v(x)\,dm(x)$$
for any function $\varphi\in C_1^{\,0}(U).$ Then we say that the
function $v$ is a weak derivative of the first order of the function
$u$ with respect to $x_i$ and denote it by $\frac{\partial
u}{\partial x_i}(x):= v.$

A function $u\in W_{\rm loc}^{1,1}(U)$ if $u$ has weak derivatives
of the first order with respect to each of the variables in $U,$
which are locally integrable in $U.$

\medskip
A mapping $f:U\rightarrow {\Bbb R}^n$ belongs to the Sobolev class
$W_{\rm loc}^{1,1}(U),$ write $f \in W^{1,1}_{\rm loc}(U),$ if all
coordinate functions of $f=(f_1,\ldots, f_n)$ have weak partial
derivatives of the first order, which are locally integrable in $U.$
We write $f\in W^{1, p}_{\rm loc}(U),$ $p\geqslant 1,$ if all
coordinate functions of $f=(f_1,\ldots, f_n) $ have weak partial
derivatives of the first order, which are locally integrable in $U$
to the degree $p.$

\medskip
Let $\varphi:{\Bbb R}^+\rightarrow{\Bbb R}^+$ be a measurable
function. The {\it Orlicz-Sobolev class} $W^{1,\varphi}_{\rm
loc}(D)$ is the class of all locally integrable functions $f$ with
the first weak derivatives whose gradient $\nabla f$ belongs locally
in $D$ to the Orlicz class. Note that by the definition
$W^{1,\varphi}_{\rm loc}\subset W^{1,1}_{\rm loc}$. For the case
when $\varphi(t)=t^p$, $p\geqslant 1,$ we write as usual $f\in
W^{1,p}_{\rm loc}.$

\medskip
Later on, we also write $f\in W^{1,\varphi}_{\rm loc}(D)$ for a
locally integrable vector-function $f=(f_1,\ldots,f_m)$ of $n$ real
variables $x_1,\ldots,x_n$ if $f_i\in W^{1,1}_{\rm loc}(D)$ and
$$
\int\limits_{G}\varphi\left(|\nabla f(x)|\right)\,dm(x)<\infty
$$
for any domain $G\subset D$ such that $\overline{G}\subset D,$ where
$|\nabla f(x)|=\sqrt{\sum\limits_{i,j}\left(\frac{\partial
f_i}{\partial x_j}\right)^2}.$

\medskip
Recall that a mapping $f:D\rightarrow {\Bbb R}^n$ is called {\it
discrete} if the pre-image $\{f^{-1}\left(y\right)\}$ of each point
$y\,\in\,{\Bbb R}^n$ consists of isolated points, and {\it is open}
if the image of any open set $U\subset D$ is an open set in ${\Bbb
R}^n.$

\medskip
The boundary of $D$ is called {\it weakly flat} at the point $x_0\in
\partial D,$ if for every $P>0$ and for any neighborhood $U$
of the point $x_0$ there is a neighborhood $V\subset U$ of the same
point such that $M(\Gamma(E, F, D))>P$ for any continua $E, F\subset
D$ such that $E\cap\partial U\ne\varnothing\ne E\cap\partial V$ and
$F\cap\partial U\ne\varnothing\ne F\cap\partial V.$ The boundary of
$D$ is called weakly flat if the corresponding property holds at any
point of the boundary of a domain $D.$

\medskip
Given a mapping $f:D\rightarrow {\Bbb R}^n$, we denote
\begin{equation}\label{eq1_A_4} C(f, x):=\{y\in \overline{{\Bbb
R}^n}:\exists\,x_k\in D: x_k\rightarrow x, f(x_k) \rightarrow y,
k\rightarrow\infty\}
\end{equation}
and
\begin{equation}\label{eq1_A_5} C(f, \partial
D)=\bigcup\limits_{x\in \partial D}C(f, x)\,.
\end{equation}
In what follows, ${\rm Int\,}A$ denotes the set of inner points of
the set $A\subset \overline{{\Bbb R}^n}.$ Recall that, the set
$U\subset\overline{{\Bbb R}^n}$ is neighborhood of the point $z_0,$
if $z_0\in {\rm Int\,}A.$

\medskip
We say that a function ${\varphi}:D\rightarrow{\Bbb R}$ has a {\it
finite mean oscillation} at a point $x_0\in D,$ write $\varphi\in
FMO(x_0),$ if
\begin{equation}\label{eq29*!}
\limsup\limits_{\varepsilon\rightarrow
0}\frac{1}{\Omega_n\varepsilon^n}\int\limits_{B( x_0,\,\varepsilon)}
|{\varphi}(x)-\overline{{\varphi}}_{\varepsilon}|\ dm(x)<\infty\,,
\end{equation}
where
$$\overline{{\varphi}}_{\varepsilon}=\frac{1}
{\Omega_n\varepsilon^n}\int\limits_{B(x_0,\,\varepsilon)}
{\varphi}(x) \,dm(x)\,.$$
Let $Q:{\Bbb R}^n\rightarrow [0,\infty]$ be a Lebesgue measurable
function. We set
$$Q^{\,\prime}(x)=\left\{
\begin{array}{rr}
Q(x), &   Q(x)\geqslant 1\,, \\
1,  &  Q(x)<1\,.
\end{array}
\right.$$ Denote by $q^{\,\prime}_{x_0}$ the mean value of
$Q^{\,\prime}(x)$ over the sphere $|x-x_0|=r$, that means,
\begin{equation}\label{eq32*B}
q^{\,\prime}_{x_0}(r):=\frac{1}{\omega_{n-1}r^{n-1}}
\int\limits_{|x-x_0|=r}Q^{\,\prime}(x)\,d{\mathcal H}^{n-1}\,.
\end{equation}
Using the inversion $\psi(x)=\frac{x}{|x|^2},$ we may define $FMO$
for $x_0=\infty.$

\medskip
We say that the boundary $\partial D$ of a~domain $D$ in ${\Bbb
R}^n,$ $n\geqslant 2,$ is {\it strongly accessible at a~point
$x_0\in
\partial D$ with respect to the $p$-modulus} if for each neighborhood
$U$ of $x_0$ there
exist a~compact set $E\subset D$, a~neighborhood $V\subset U$ of
$x_0$ and $\delta>0$ such that
\begin{equation}
\label{eq1.3_a} M_p(\Gamma(E,F, D))\geqslant \delta
\end{equation}
for each continuum $F$ in~$D$ that intersects $\partial U$ and
$\partial V$. When $p=n$, we will usually drop the prefix in the
``$p$-modulus'' when speaking about~(\ref{eq1.3_a}).

\medskip
Assume that, a mapping $f$ has partial derivatives almost everywhere
in $D.$ In this case, we set
$$l\left(f^{\,\prime}(x)\right)\,=\,\min\limits_{h\in {\Bbb R}^n
\backslash \{0\}} \frac {|f^{\,\prime}(x)h|}{|h|}\,,$$
\begin{equation}\label{eq5_a}
\Vert f^{\,\prime}(x)\Vert\,=\,\max\limits_{h\in {\Bbb R}^n
\backslash \{0\}} \frac {|f^{\,\prime}(x)h|}{|h|}\,,
\end{equation}
$$J(x,f)=\det
f^{\,\prime}(x)\,.$$
and define for any $x\in D$ and $\alpha\geqslant 1$
\begin{equation}\label{eq0.1.1A}
K_{I, \alpha}(x,f)\quad =\quad\left\{
\begin{array}{rr}
\frac{|J(x,f)|}{{l\left(f^{\,\prime}(x)\right)}^{\,\alpha}}, & J(x,f)\ne 0,\\
1,  &  f^{\,\prime}(x)=0, \\
\infty, & {\rm otherwise}
\end{array}
\right.\,,\end{equation}
The following statement holds, cf.~\cite[Lemma~5.20,
Corollary~5.23]{MRSY$_1$}, \cite[Lemma~6.1, Theorem~6.1]{RS$_1$},
\cite[Lemma~5, Theorem~3]{Sm} and \cite[Theorem~1.1]{DS}.

\medskip
\begin{theorem}\label{th3}
{\it\, Let $\alpha\geqslant 1,$ let $D$ and $D^{\,\prime}$ be
bounded domains in ${\Bbb R}^n,$ $n\geqslant 3,$ let $b\in \partial
D,$ let $Q:D\rightarrow [0, \infty]$ be a Lebesgue measurable
function and let $\varphi:[0,\infty)\rightarrow[0,\infty)$ be an
increasing function. Let $f:D\rightarrow D^{\,\prime}$ be a bounded
open discrete mapping of the class $W^{1,\varphi}_{\rm loc}(D),$
$f(D)=D^{\,\prime}.$ In addition, assume that $C(f, \partial
D)\subset E_*$ for some closed (with respect to $D^{\,\prime}$) set
$E_*\subset D^{\,\prime}$ and $f^{\,-1}(E_*)=E$ for some closed
(with respect to $D$) set $E\subset D,$ besides that,

\medskip
1) the set $E$
is nowhere dense in $D$ and $D$ is finitely connected on $E\cup
\partial D,$ i.e., for any $z_0\in E\cup \partial D$ and any neighborhood $\widetilde{U}$
of $z_0$ there is a neighborhood $\widetilde{V}\subset
\widetilde{U}$ of $z_0$ such that $(D\cap \widetilde{V})\setminus E$
consists of finite number of components,

\medskip
2) for any neighborhood $U$ of $b$ there is a neighborhood $V\subset
U$ of $b$ such that:

\medskip
2a) $V\cap D$ is connected,

\medskip
2b) $(V\cap D)\setminus E$ consists at most of $m$ components,
$1\leqslant m<\infty,$

\medskip
3) all components of the set $D^{\,\prime}\setminus E_*$ have a
strongly accessible boundary with respect to $\alpha$-modulus,

\medskip
4) the function $\varphi$ satisfies the following Calderon condition
\begin{equation}\label{eq1_A_10}
\int\limits_{1}^{\infty}\left(\frac{t}{\varphi(t)}\right)^
{\frac{1}{n-2}}\,dt<\infty\,.
\end{equation}

Suppose that $K_{I, \alpha}(x, f)\leqslant Q(x)$ a.e. and at least
one of the following conditions is satisfied:

\medskip
5.1) the function $Q$ has a finite mean oscillation at a point $b;$

\medskip
5.2) $q_{b}(r)\,=\,O\left(\left[\log{\frac1r}\right]^{n-1}\right)$
as $r\rightarrow 0;$

\medskip 5.3) the condition
\begin{equation}\label{eq6}
\int\limits_{0}^{\delta(b)}\frac{dt}{t^{\frac{n-1}{\alpha-1}}
q_{b}^{\,\prime\,\frac{1}{\alpha-1}}(t)}=\infty
\end{equation}
holds for some $\delta(b)>0.$ Then $f$ has a continuous extension to
$b.$

\medskip
If the above is true for any point $b\in\partial D,$ the mapping $f$
has a continuous extension
$\overline{f}:\overline{D}\rightarrow\overline{D^{\,\prime}},$
moreover, $\overline{f}(\overline{D})=\overline{D^{\,\prime}}.$ }
\end{theorem}

\medskip
For Sobolev classes on the plane, Theorem~\ref{th3} has a somewhat
simpler form.

\medskip
\begin{theorem}\label{th1}
{\it\, Let $\alpha\geqslant 1,$ let $D$ and $D^{\,\prime}$ be
bounded domains in ${\Bbb R}^2,$ let $b\in \partial D,$ let
$Q:D\rightarrow [0, \infty]$ be a locally integrable function. Let
$f:D\rightarrow D^{\,\prime}$ be a bounded open discrete mapping of
the class $W^{1,1}_{\rm loc}(D),$ $f(D)=D^{\,\prime}.$ In addition,
assume that $C(f, \partial D)\subset E_*$ for some closed (with
respect to $D^{\,\prime}$) set $E_*\subset D^{\,\prime}$ and
$f^{\,-1}(E_*)=E$ for some closed (with respect to $D$) set
$E\subset \overline{D},$ besides that, the conditions 1)--3) from
Theorem~\ref{th3} hold. Suppose also that $K_{I, \alpha}(x,
f)\leqslant Q(x)$ a.e. and at least one of the following conditions
is satisfied: 1) the function $Q$ has a finite mean oscillation at a
point $b;$ 2)
$q_{b}(r)\,=\,O\left(\left[\log{\frac1r}\right]\right)$ as
$r\rightarrow 0;$ 3) the condition~(\ref{eq6}) holds for some
$\delta(b)>0.$ Then $f$ has a continuous extension to $b.$

If the above is true for any point $b\in\partial D,$ the mapping $f$
has a continuous extension
$\overline{f}:\overline{D}\rightarrow\overline{D^{\,\prime}},$
moreover, $\overline{f}(\overline{D})=\overline{D^{\,\prime}}.$ }
\end{theorem}

\section{Preliminaries}

A Borel function $\rho:{\Bbb R}^n\,\rightarrow [0,\infty] $ is
called {\it admissible} for the family $\Gamma$ of paths $\gamma$ in
${\Bbb R}^n,$ if the relation
\begin{equation*}\label{eq1.4}
\int\limits_{\gamma}\rho (x)\, |dx|\geqslant 1
\end{equation*}
holds for all (locally rectifiable) paths $ \gamma \in \Gamma.$ In
this case, we write: $\rho \in {\rm adm} \,\Gamma .$ Let $p\geqslant
1,$ then {\it $p$-modulus} of $\Gamma $ is defined by the equality
\begin{equation*}\label{eq1.3gl0}
M_p(\Gamma)=\inf\limits_{\rho \in \,{\rm adm}\,\Gamma}
\int\limits_{{\Bbb R}^n} \rho^p (x)\,dm(x)\,.
\end{equation*}
Let $x_0\in {\Bbb R}^n,$ $0<r_1<r_2<\infty,$
\begin{equation}\label{eq1ED}
S(x_0,r) = \{ x\,\in\,{\Bbb R}^n : |x-x_0|=r\}\,, \quad B(x_0, r)=\{
x\,\in\,{\Bbb R}^n : |x-x_0|<r\}\end{equation}
and
\begin{equation}\label{eq1**}
A=A(x_0, r_1,r_2)=\left\{ x\,\in\,{\Bbb R}^n:
r_1<|x-x_0|<r_2\right\}\,.\end{equation}
The most important tool of the paper is the connection between
Sobolev (Orlicz-Sobolev) classes and lower (ring) $Q$-mappings. The
theoretical part of this connection is mainly established
in~\cite{Sev$_2$}.

\medskip
Let $\omega$ be an open set in $\overline{{\Bbb R}^k}:={\Bbb
R}^k\cup\{\infty\}$, $k=1,\ldots\,,n-1$. A (continuous) mapping
$S:\omega\rightarrow{\Bbb R}^n$ is called a $k$-dimensional surface
$S$ in ${\Bbb R}^n$. Sometimes we call the image
$S(\omega)\subseteq{\Bbb R}^n$ the surface $S$, too. The number of
preimages
\begin{equation} \label{eq8.2.3}N(S,y)\ =\
\mathrm{card}\,\:S^{-1}(y)=\mathrm{card}\:\{x\in\omega:\
S(x)=y\}\end{equation}
is said to be a {\it multiplicity function} of the surface $S$ at a
point $y\in{\Bbb R}^n$. In other words, $N(S,y)$ denotes the
multiplicity of covering of the point $y$ by the surface $S$. It is
known that the multiplicity function is lower semi continuous, i.e.,
$$N(S,y) \geqslant \liminf_{m\rightarrow\infty}: N(S,y_m)$$
for every sequence $y_m\in{\Bbb R}^n$, $m=1,2,\ldots\,$, such that
$y_m\rightarrow y\in{\Bbb R}^n$ as $m\rightarrow\infty$, see, e.g.,
\cite{RR}, p. 160. Thus, the function $N(S,y)$ is Borel measurable
and hence measurable with respect to every Hausdorff measure ${\cal
H}^k$ (see, e.g., \cite{Sa}, p.~52).

\medskip
If $\rho:{\Bbb R}^n\rightarrow[0,\infty]$ is a Borel function, then
its {\it integral over a $k$-dimensional surface} $S$ in ${\Bbb
R}^n$, $n\geqslant 2,$ is defined by the equality
\begin{equation}\label{eq8.2.5} \int\limits_S \rho\ d{\cal {A}}\
:=\ \int\limits_{{\Bbb R}^n}\rho(y)\:N(S,y)\ d{\cal H}^ky\
.\end{equation} Given a family ${\cal S}$ of such $k$-dimensional
surfaces $S$ in ${\Bbb R}^n$, a Borel function $\ \rho:{\Bbb
R}^n\rightarrow[0,\infty]$ is called {\it admissible} for $\cal S$,
abbr. $\rho\in\mathrm{adm}\,\cal S$, if
\begin{equation}\label{eq8.2.6A}
\int\limits_S\rho^k\ d{\cal{A}}\ \geqslant\ 1\end{equation} for
every $S\in\cal S$. Given $p\in[k,\infty)$, the {\it $p$-modulus} of
$\cal S$ is the quantity
\begin{equation}\label{M} M_p({\cal S})\ =\
\inf_{\rho\in\mathrm{adm}\,\cal S}\int\limits_{{\Bbb R}^n}\rho^p(x)\
dm(x)\,.\end{equation}
The next class of mappings is a generalization of quasiconformal
mappings in the sense of Gehring's ring definition (see \cite{Ge};
cf. \cite[Chapter~9]{MRSY$_2$}). Let $D$ and $D^{\,\prime}$ be
domains in ${\Bbb R}^n,$ $n\geqslant 2$. Suppose that
$x_0\in\overline {D}\setminus\{\infty\}$ and $Q\colon
D\rightarrow(0,\infty)$ is a Lebesgue measurable function. A mapping
$f:D\rightarrow D^{\,\prime}$ is called a {\it lower $Q$-mapping at
a point $x_0$ relative to $p$-modulus} if
\begin{equation}\label{eq1A}
M_p(f(\Sigma_{\varepsilon}))\geqslant \inf\limits_{\rho\in{\rm
ext}_p\,{\rm adm}\Sigma_{\varepsilon}}\int\limits_{D\cap A(x_0,
\varepsilon, r_0)}\frac{\rho^p(x)}{Q(x)}\,dm(x)
\end{equation}
for every spherical ring $A(x_0, \varepsilon, r_0)=\{x\in {\Bbb
R}^n\,:\, \varepsilon<|x-x_0|<r_0\}$, $r_0\in(0,d_0)$,
$d_0=\sup\limits_{x\in D}|x-x_0|$, where $\Sigma_{\varepsilon}$ is
the family of all intersections of the spheres $S(x_0, r)$ with the
domain $D$, $r\in (\varepsilon, r_0)$. If $p=n$, we say that $f$ is
a lower $Q$-mapping at $x_0$. We say that $f$ is a lower $Q$-mapping
relative to $p$-modulus in $A\subset {\Bbb R}^n$ if (\ref{eq1A}) is
true for all $x_0\in A$.

\medskip
Let
\begin{equation}\label{eq1_A_2}
N(y, f, A)\,=\,{\rm card}\,\left\{x\in A: f(x)=y\right\}\,,
\end{equation}
\begin{equation}\label{eq1_A_3}
N(f, A)\,=\,\sup\limits_{y\in{\Bbb R}^n}\,N(y, f, A)\,.
\end{equation}
The following lemma holds, see e.g. \cite[Lemma~5.1]{Sev$_1$}.
\begin{lemma}\label{thOS4.1} {\it\,Let $D$ be a domain in ${\Bbb R}^n$,
$n\geqslant  2$, and let $\varphi\colon (0,\infty)\rightarrow
(0,\infty)$ be a monotone nondecreasing function satisfying
(\ref{eq1_A_10}). If $p>n-1$, then every open discrete mapping
$f\colon D\rightarrow {\Bbb R}^n$ of class $W^{1,\varphi}_{{\rm
loc}}$ with finite distortion and such that $N(f, D)<\infty$ is a
lower $Q$-mapping relative to the $p$-modulus at every point
$x_0\in\overline {D}$ for
$$
Q(x)=N(f, D)\cdot K^{\frac{p-n+1}{n-1}}_{I, \alpha}(x, f),
$$
$\alpha:=\frac{p}{p-n+1}$, where the inner dilation $K_{I,\alpha}(x,
f)$ for~$f$ at~$x$ of order~$\alpha$ is defined by \eqref{eq0.1.1A},
and the multiplicity $N(f, D)$ is defined by the
relation~(\ref{eq1_A_3}).}
\end{lemma}

\medskip
The statement similar to Lemma~\ref{thOS4.1} holds for $n=2,$ but
for Sobolev classes (see e.g. \cite[Theorem~4]{Sev$_3$}).

\medskip
\begin{lemma}\label{lem3} {\it\,Let $D$ be a domain in ${\Bbb R}^2$,
$n\geqslant  2$, and let $p>1.$ Then every open discrete mapping
$f\colon D\rightarrow {\Bbb R}^2$ of class $W^{1,1}_{{\rm loc}}$
with finite distortion and such that $N(f, D)<\infty$ is a lower
$Q$-mapping relative to the $p$-modulus at every point
$x_0\in\overline {D}$ for
$$
Q(x)=N(f, D)\cdot K^{p-1}_{I, \alpha}(x, f),
$$
$\alpha:=\frac{p}{p-1}$, where the inner dilation $K_{I,\alpha}(x,
f)$ for~$f$ at~$x$ of order~$\alpha$ is defined by \eqref{eq0.1.1A},
and the multiplicity $N(f, D)$ is defined by the
relation~(\ref{eq1_A_3}).}
\end{lemma}

\medskip
Given sets $E,$ $F\subset\overline{{\Bbb R}^n}$ and a domain
$D\subset {\Bbb R}^n$ we denote by $\Gamma(E,F,D)$ a family of all
paths $\gamma:[a,b]\rightarrow \overline{{\Bbb R}^n}$ such that
$\gamma(a)\in E,\gamma(b)\in\,F $ and $\gamma(t)\in D$ for $t \in
(a, b).$ Let $S_i=S(x_0, r_i),$ $i=1,2,$ where spheres $S(x_0, r_i)$
centered at $x_0$ of the radius $r_i$ are defined in~(\ref{eq1ED}).
The following statement holds, see~\cite{Sev$_2$}.

\medskip
\begin{lemma}\label{lem4}
{\,\it Let $x_0\in \partial D,$ let $f:D\rightarrow {\Bbb R}^n$ be a
bounded, open, discrete, and closed lower $Q$-mapping with respect
to $p$-modulus in a domain $D\subset{\Bbb R}^n,$ $Q\in
L_{loc}^{\frac{n-1}{p-n+1}}({\Bbb R}^n),$ $n-1<p,$ and
$\alpha:=\frac{p}{p-n+1}.$ Then, for every
$\varepsilon_0<d_0:=\sup\limits_{x\in D}|x-x_0|$ and every compact
set $C_2\subset D\setminus B(x_0, \varepsilon_0)$ there exists
$\varepsilon_1,$ $0<\varepsilon_1<\varepsilon_0,$ such that, for
each $\varepsilon\in (0, \varepsilon_1)$ and each compact
$C_1\subset \overline{B(x_0, \varepsilon)}\cap D$ the inequality
\begin{equation}\label{eq3A_2}
M_{\alpha}(f(\Gamma(C_1, C_2, D)))\leqslant \int\limits_{A(x_0,
\varepsilon, \varepsilon_1)}Q^{\frac{n-1}{p-n+1}}(x)
\eta^{\alpha}(|x-x_0|)\,dm(x)
\end{equation}
holds, where $A(x_0, \varepsilon, \varepsilon_1)=\{x\in {\Bbb R}^n:
\varepsilon<|x-x_0|<\varepsilon_1\}$ and $\eta: (\varepsilon,
\varepsilon_1)\rightarrow [0,\infty]$ is an arbitrary Lebesgue
measurable function such that
\begin{equation}\label{eq6B}
\int\limits_{\varepsilon}^{\varepsilon_1}\eta(r)\,dr=1\,.
\end{equation}
}
\end{lemma}
Let $D\subset {\Bbb R}^n,$ let $f:D\rightarrow {\Bbb R}^n$ be a
discrete open mapping, $\beta: [a,\,b)\rightarrow {\Bbb R}^n$ be a
path, and $x\in\,f^{\,-1}(\beta(a)).$ A path $\alpha:
[a,\,c)\rightarrow D$ is called a {\it maximal $f$-lifting} of
$\beta$ starting at $x,$ if $(1)\quad \alpha(a)=x\,;$ $(2)\quad
f\circ\alpha=\beta|_{[a,\,c)};$ $(3)$\quad for
$c<c^{\prime}\leqslant b,$ there is no a path $\alpha^{\prime}:
[a,\,c^{\prime})\rightarrow D$ such that
$\alpha=\alpha^{\prime}|_{[a,\,c)}$ and $f\circ
\alpha^{\,\prime}=\beta|_{[a,\,c^{\prime})}.$ Here and below we say
that a path $\beta:[a, b)\rightarrow \overline{{\Bbb R}^n}$
converges to the set $C\subset \overline{{\Bbb R}^n}$ as
$t\rightarrow b,$ if $h(\beta(t), C)=\sup\limits_{x\in C}h(\beta(t),
C)\rightarrow 0$ at $t\rightarrow b.$ The following is true
(see~\cite[Lemma~3.12]{MRV}).

\medskip
\begin{proposition}\label{pr3}
{\it\, Let $f:D\rightarrow {\Bbb R}^n,$ $n\geqslant 2,$ be an open
discrete mapping, let $x_0\in D,$ and let $\beta: [a,\,b)\rightarrow
{\Bbb R}^n$ be a path such that $\beta(a)=f(x_0)$ and such that
either $\lim\limits_{t\rightarrow b}\beta(t)$ exists, or
$\beta(t)\rightarrow \partial f(D)$ as $t\rightarrow b.$ Then
$\beta$ has a maximal $f$-lifting $\alpha: [a,\,c)\rightarrow D$
starting at $x_0.$ If $\alpha(t)\rightarrow x_1\in D$ as
$t\rightarrow c,$ then $c=b$ and $f(x_1)=\lim\limits_{t\rightarrow
b}\beta(t).$ Otherwise $\alpha(t)\rightarrow \partial D$ as
$t\rightarrow c.$}
\end{proposition}

\medskip
The following statement is proved in~\cite[Lemma~2.1]{DS}.

\medskip
\begin{lemma}\label{lem2}
{\it\, Let $D$ be a domain in ${\Bbb R}^n,$ $n\geqslant 2,$ and let
$x_0\in \partial D.$ Assume that $E$ is closed and nowhere dense in
$D,$ and $D$ is finitely connected on $E,$ i.e., for any $z_0\in E$
and any neighborhood $\widetilde{U}$ of $z_0$ there is a
neighborhood $\widetilde{V}\subset \widetilde{U}$ of $z_0$ such that
$(D\cap \widetilde{V})\setminus E$ consists of finite number of
components.

In addition, assume that the following condition holds: for any
neighborhood $U$ of $x_0$ there is a neighborhood $V\subset U$ of
$x_0$ such that:

\medskip
a) $V\cap D$ is connected,

\medskip
b) $(V\cap D)\setminus E$ consists at most of $m$ components,
$1\leqslant m<\infty.$

Let $x_k, y_k\in D\setminus E,$ $k=1,2,\ldots ,$ be a sequences
converging to $x_0$ as $k\rightarrow\infty.$ Then there are
subsequences $x_{k_l}$ and $y_{k_l},$ $l=1,2,\ldots ,$ belonging to
some sequence of neighborhoods $V_l,$ $l=1,2,\ldots ,$ of the point
$x_0$ such that $V_l\subset B(x_0, 2^{\,-l}),$ $l=1,2,\ldots ,$ and,
in addition, any pair $x_{k_l}$ and $y_{k_l}$ may be joined by a
path $\gamma_l$ in $V_l\cap D,$ where $\gamma_l$ contains at most
$m-1$ points in~$E.$ }
\end{lemma}

The following statement may be found in~\cite[Lemmas~1.3 and
1.4]{Sev$_1$}.

\medskip
\begin{proposition}\label{pr6}
{\it\, Let $Q:{\Bbb R}^n\rightarrow [0,\infty],$ $n\geqslant 2,$
$n-1<p\leqslant n,$ be a Lebesgue measurable function and let
$x_0\in {\Bbb R}^n.$ Assume that either of the following conditions
holds

\noindent (a) $Q\in FMO(x_0),$

\noindent (b)
$q_{x_0}(r)\,=\,O\left(\left[\log{\frac1r}\right]^{n-1}\right)$ as
$r\rightarrow 0,$

\noindent (c) for some small $\delta_0=\delta_0(x_0)>0$ we have the
relations
$$
\int\limits_{\delta}^{\delta_0}\frac{dt}{t^{\frac{n-1}{p-1}}q_{x_0}^{\frac{1}{p-1}}(t)}<\infty,\qquad
0<\delta<\delta_0,
$$
and
$$
\int\limits_{0}^{\delta_0}\frac{dt}{t^{\frac{n-1}{p-1}}q_{x_0}^{\frac{1}{p-1}}(t)}=\infty\,.
$$
Then  there exist a number $\varepsilon_0\in(0,1)$ and a function
$\psi:(0, \varepsilon_0)\rightarrow [0, \infty)$ such that the
relation
$$
\int\limits_{\varepsilon<|x-x_0|<\varepsilon_0}Q(x)\cdot\psi^p(|x-x_0|)
\ dm(x)=o(I^p(\varepsilon, \varepsilon_0))
$$
holds as $\varepsilon\rightarrow 0,$
where
$$
0<I(\varepsilon, \varepsilon_0)
=\int\limits_{\varepsilon}^{\varepsilon_0}\psi(t)\,dt < \infty
\qquad \forall\quad\varepsilon \in(0, \varepsilon_1)
$$
for some $0<\varepsilon_1<\varepsilon_0.$ }
\end{proposition}

\section{Main Lemma}

Various versions of the lemma given below have been given in
publications \cite[Lemma~2.1]{IR}, \cite[Lemma~5.16]{MRSY$_1$},
\cite[Theorem~4.6, Theorem~13.1]{MRSY$_2$},
\cite[Theorem~5.1]{RS$_1$}, \cite{Sev$_2$} and \cite[Theorem~1]{Sm},
cf.~~\cite[Theorem~17.15]{Va}, \cite[Section~3]{Na$_1$},
\cite[Theorem~4.10.II]{Vu} and~\cite[Theorem~4.2]{Sr},
cf.~\cite{DS}.

\medskip
\begin{lemma}\label{lem1} {\it
Assume that, under conditions of Theorem~\ref{th3}, instead
of~5.1)--5.3) the following condition holds: there is
$\varepsilon_0=\varepsilon_0(b)>0$ and some positive measurable
function $\psi:(0, \varepsilon_0)\rightarrow (0,\infty)$ such that
\begin{equation}\label{eq7***}
0<I(\varepsilon,
\varepsilon_0)=\int\limits_{\varepsilon}^{\varepsilon_0}\psi(t)\,dt
< \infty
\end{equation}
for sufficiently small $\varepsilon\in(0, \varepsilon_0)$ and, in
addition,
\begin{equation}\label{eq5***}
\int\limits_{A(b, \varepsilon, \varepsilon_0)}
Q(x)\cdot\psi^{\,\alpha}(|x-b|)
 \ dm(x) =o(I^{\alpha}(\varepsilon, \varepsilon_0))\,,
\end{equation}
where $A:=A(b, \varepsilon, \varepsilon_0)$ is defined
in~(\ref{eq1**}). Then $f$ has a continuous extension to $b.$}
\end{lemma}

\medskip
\begin{proof}
Firstly let us to prove that, under conditions of the lemma,
$D\setminus E$ consists of finite number of components $D_1,\ldots,
D_s,$ $1\leqslant s<\infty.$ Indeed, assume the contrary, namely,
that $D\setminus E$ consists of infinite number of components $D_1,
D_2,\ldots .$ Let $z_i\in D_i,$ $i=1,2,\ldots .$ Since $D$ is
bounded, there is a subsequence $z_{i_k},$ $k_1,2,\ldots ,$ which
converge to some point $z_0\in \overline{D}.$ On the other hand,
there is a neighborhood $V$ of the point $z_0$ such that $(V\cap
D)\setminus E$ consists of $m$ components. Thus, there is at least
one such a component $K$ intersecting infinitely many components
$D_{i_1}, D_{i_2},\ldots .$ However, by the assumption there is a
neighborhood $V$ of $z_0$ such that $(V\cap D)\setminus E$ consists
of finite number of components $V_1,\ldots , V_m.$ Now, there is
$m_0\in[1, m]$ such that $V_{m_0}$ intersects infinitely many $D_i,$
contradiction. Thus, $D_1,\ldots, D_s,$ $1\leqslant s<\infty,$ as
required.

\medskip
Now, let us to prove that $f$ is closed in each $D_i,$
$i=1,2,\ldots, $ i.e. $f(S)$ is closed in $f(D_i)$ whenever $S$ is
closed in $D_i.$ Since $f$ is open and discrete, it is sufficient to
prove that $f$ is boundary preserving, that is $C(f,
\partial D_i)\subset \partial f(D_i)$ (see, e.g.,
\cite[Theorem~3.3]{Vu}). Indeed, let $z_k\in D_i,$ $k=1,2,\ldots ,$
let $z_k\rightarrow z_0\in\partial D_i$ and let $f(z_k)\rightarrow
w_0$ as $k\rightarrow\infty.$ We need to prove that $w_0\in
\partial f(D_i).$ Assume the contrary, i.e. $w_0\in f(D_i).$
There are two cases: $z_0\in \partial D$ or $z_0\in D.$ 1) In the
first case, when $z_0\in \partial D,$ we obtain that $w_0\in C(f,
\partial D)\subset E_*.$ But since $w_0\in f(D_i),$ there is
$\zeta_i\in D_i$ such that $f(\zeta_i)=w_0.$ Now, since $w_0\in
E_*,$ we obtain that $\zeta_i\in f^{\,-1}(E_*)=E$ and, consequently,
$\zeta_i\not\in D_i,$ contradiction. 2) Let us consider the second
case, when $z_0\in D.$ Now, by the definition of $D_i,$ we have that
$z_0\in E.$ Now, $w_0\in f(E)\subset E_*.$ But since $w_0\in
f(D_i),$ there is $\zeta_i\in D_i$ such that $f(\zeta_i)=w_0.$ Since
$w_0\in E_*,$ we have that $\zeta_i\in f^{\,-1}(E_*)=E$ and,
consequently, $\zeta_i\not\in D_i,$ contradiction.

\medskip
Let $K$ be a component of $D^{\,\prime}\setminus E_*$ consisting
$f(D_i).$ Observe that $f(D_i)=K.$ Indeed, $f(D_i)\subset K$ by the
definition. Let us to prove that $K\subset f(D_i).$ Let us prove
this inclusion by the contradiction, i.e., let $a_0\in K\setminus
f(D_i).$ Chose $b_0\in f(D_i)$ and join the points $b_0$ and $a_0$
by a path $\beta:[0, 1]\rightarrow K.$ Let $\alpha:[0, c)\rightarrow
D$ be a maximal $f$-lifting of $\beta$ starting at
$c_0:=f^{\,-1}(b_0)\cap D_i$ (this lifting exists by
Proposition~\ref{pr3}). By the same proposition either one of two
cases are possible: $\alpha(t)\rightarrow x_1$ as $t\rightarrow c,$
or $\alpha(t)\rightarrow \partial D_i$ as $t\rightarrow c.$ In the
first case, by Proposition~\ref{pr3} we obtain that $c=1$ and
$f(\beta(1))=f(x_1)=a_0$ which contradict with the choice of $a_0.$
In the second case, when $\alpha(t)\rightarrow \partial D_i$ as
$t\rightarrow c,$ we obtain that $f(\beta(c))\in C(f, \partial
D_i)\subset E_*.$ The latter contradicts the definition of $\beta$
because $\beta$ does not contain itself points in $E_*.$

\medskip
Observe that, $D^{\,\prime}\setminus E_*$ consists of finite number
of components. Otherwise, $D^{\,\prime}\setminus
E_*=\bigcup\limits_{i=1}^{\infty}K^{\,\prime}_i,$ $\varnothing \ne
K^{\,\prime}_i\ne K^{\,\prime}_j\ne \varnothing$ for $i\ne j.$  Let
$K_i$ be some component of the set $f^{\,-1}(K^{\,\prime}_i).$ By
the definition $K_i\subset D\setminus E$ and $\varnothing \ne K_i\ne
K_j\ne \varnothing$ for $i\ne j.$ The latter contradicts the proving
above for the number of components of $D\setminus E.$

\medskip
Now, we proceed directly to the proof of the main statement of the
lemma. We apply the approach used under the proof of
\cite[Lemma~3.1]{DS}. Suppose the opposite, i.e., the conclusion of
the lemma is not true. Since $D^{\,\prime}$ is bounded, there are at
least two sequences $x_k,$ $y_k\in D,$ $i=1,2,\ldots,$ such that
$x_k,$ $y_k\in D,$ $k=1,2,\ldots,$ $x_k\rightarrow b,$
$y_k\rightarrow b$ as $k\rightarrow \infty,$ and $f(x_k)\rightarrow
y,$ $f(y_k)\rightarrow y^{\,\prime}$ as $k\rightarrow \infty,$ while
$y^{\,\prime}\ne y.$
In particular,
\begin{equation}\label{eq1B}
|f(x_k)-f(y_k)|\geqslant \delta>0
\end{equation}
for some $\delta>0$ and all $k\in {\Bbb N}.$
By the assumption~2), there exists a sequence of neighborhoods
$V_k\subset B(b, 2^{\,-k}),$ $k=1,2,\ldots ,$ such that $V_k\cap D$
is connected and $(V_k\cap D)\setminus E$ consists at most of $m$
components, $1\leqslant m<\infty.$

We note that the points $x_k$ and $y_k,$ $k=1,2,\ldots, $ may be
chosen such that $x_k, y_k\not\in E.$ Indeed, since under
condition~1) the set $E$ is nowhere dense in $D,$ there exists a
sequence $x_{ki}\in D\setminus E,$ $i=1,2,\ldots ,$ such that
$x_{ki}\rightarrow x_k$ as $i\rightarrow\infty.$ Put
$\varepsilon>0.$ Due to the continuity of the mapping $f$ at the
point $x_k,$ for the number $k\in {\Bbb N}$ there is a number
$i_k\in {\Bbb N}$ such that $|f(x_{ki_k})-f(x_k)|<\frac{1}{2^k}.$
So, by the triangle inequality
$$|f(x_{ki_k})-y|\leqslant |f(x_{ki_k})-f(x_k)|+
|f(x_k)-y|\leqslant \frac{1}{2^k}+\varepsilon\,,$$
$k\geqslant k_0=k_0(\varepsilon),$ because $f(x_k)\rightarrow y$ as
$k\rightarrow\infty$ and by the choice of $x_k$ and $y.$ Therefore,
$x_k\in D$ may be replaced by $x_{ki_k}\in D\setminus E,$ as
required. We may reason similarly for the sequence $y_k.$

\medskip
Now, by Lemma~\ref{lem2} there are subsequences $x_{k_l}$ and
$y_{k_l},$ $l=1,2,\ldots ,$ belonging to some sequence of
neighborhoods $V_l,$ $l=1,2,\ldots ,$ of the point $b$ such that
${\rm diam\,}V_l\rightarrow 0$ as $l\rightarrow\infty$ and, in
addition, any pair $x_{k_l}$ and $y_{k_l}$ may be joined by a path
$\gamma_l$ in $V_l,$ where $\gamma_l$ contains at most $m-1$ points
in $E.$ Without loss of generality, we may assume that the same
sequences $x_k$ and $y_k$ satisfy properties mentioned above. Let
$\gamma_k:[0, 1]\rightarrow D,$ $\gamma_k(0)=x_k$ and
$\gamma_k(1)=y_k,$ $k=1,2,\ldots .$

\medskip
Observe that, the path $f(\gamma_k)$ contains not more than $m-1$
points in $E_*.$ In the contrary case, there are at least $m$ such
points $b_{1}=f(\gamma_k(t_1)), b_{2}=f(\gamma_k(t_2)),\ldots,
b_m=f(\gamma_k(t_m)),$ $0\leqslant t_1\leqslant t_2\leqslant
\ldots\leqslant t_m\leqslant 1.$ But now the points
$a_{1}=\gamma_k(t_1), a_{2}=\gamma_k(t_2),\ldots, a_m=\gamma_k(t_m)$
are in $f^{\,-1}(E_*)=E$ and simultaneously belong to $\gamma_k.$
This contradicts the definition of $\gamma_k.$

\medskip
Let
$$b_{1}=f(\gamma_k(t_1)), b_{2}=f(\gamma_k(t_2))\quad,\ldots,\quad
b_l=f(\gamma_k(t_1))\,,$$ $$0\leqslant t_1\leqslant t_2\leqslant
\ldots\leqslant t_l\leqslant 1,\qquad 1\leqslant l\leqslant m-1\,,$$
be points in $f(\gamma_k)\cap E_*.$ By the relation~(\ref{eq1B}) and
due to the triangle inequality,
\begin{equation}\label{eq6A}
\delta\leqslant |f(x_k)-f(y_k)|\leqslant\sum\limits_{k=1}^{l-1}
|f(\gamma_k(t_k))-f(\gamma_k(t_{k+1}))|\,.
\end{equation}
It follows from~(\ref{eq6A}) that, there is $1\leqslant
l=l(k)\leqslant m-1$ such that such that
\begin{equation}\label{eq7}
|f(\gamma_k(t_{l(k)}))-f(\gamma_k(t_{l(k)+1})|\geqslant
\delta/l\geqslant \delta/(m-1)\,.
\end{equation}
Observe that, the set $G_k:=|f(\gamma_k)|_{(t_{l(k)}, t_{l(k)+1})}|$
belongs to $D^{\,\prime}\setminus E_*,$ because it does not contain
any point in $E_*.$

Since $D^{\,\prime}\setminus E_*$ consists only of finite
components, there exists at least one a component of
$D^{\,\prime}\setminus E_*,$ containing infinitely many components
of $G_k.$ Without loss of generality, going to a subsequence, if
need, we may assume that all $G_k$ belong to one component $K$ of
$D^{\,\prime}\setminus E_*.$

\medskip
Due to the compactness of $\overline{{\Bbb R}^n},$ we may assume
that the sequence $w_k:=f(\gamma_k(t_{l(k)})),$ $k=1,2,\ldots, $
converges to some a point $w_0\in \overline{D^{\,\prime}}.$ Let us
to show that $w_0\in C(f, \partial D)\subset E_*.$ Indeed, there are
two cases: either $w_k=f(\gamma_k(t_{l(k)}))\in E_*$ for infinitely
many $k,$ or $w_k\not\in E_*$ for infinitely many $k\in {\Bbb N}.$
In the first case, the inclusion $w_0\in E_*$ is obvious, because
$E_*$ is closed. In the second case, obviously, $w_k=f(x_k),$ but
this sequence converges to $y\in C(f,
\partial D)\subset E_*$ by the assumption.

By the assumption, each component of the set $D^{\,\prime}\setminus
E_*$ has a strongly accessible boundary with respect to
$\alpha$-modulus. In this case, for any neighborhood $U$ of the
point $w_0\in
\partial K$ there is a compact set $C_0^{\,\prime}\subset
D^{\,\prime},$ a neighborhood $V$ of a point $w_0,$ $V\subset U,$
and a number $P>0$ such that
\begin{equation}\label{eq1}
M_{\alpha}(\Gamma(C_0^{\,\prime}, F, K))\geqslant P
>0
\end{equation}
for any continua $F,$ intersecting $\partial U$ and $\partial V.$
Choose a neighborhood $U$ of $w_0$ with $d(U)< \delta/2(m-1),$ where
$\delta$ is from~(\ref{eq1B}). Let $C_0^{\,\prime}$ and $V$ be a
compact set and a neighborhood corresponding to $w_0.$

Observe that, $G_k$ contains some a continuum $\widetilde{G}_k$ in
$K$ intersecting $\partial U$ and $\partial V$ for sufficiently
large $k.$ Indeed, by the construction of $G_k,$ there is a sequence
of points $a_{s, k}:=f(\gamma_k(p_s))\rightarrow
w_k:=f(\gamma_k(t_{l(k)}))$ as $s\rightarrow \infty$ and $b_{s,
k}=f(\gamma_k(q_s))\rightarrow f(\gamma_k(t_{l(k)+1}))$ as
$s\rightarrow \infty,$ where $t_{l(k)})<
p_s<q_s<\gamma_k(t_{l(k)+1}))$ and $a_{s, k}, b_{s, k}\in K.$
By the triangle inequality and by~(\ref{eq7})
$$
|a_{s, k}-b_{s, k}|\geqslant |b_{s, k}-w_k|-|w_k -a_{s, k}|
$$
\begin{equation}\label{eq11}
\geqslant|f(\gamma_k(t_{l(k)+1}))-w_k|-|f(\gamma_k(t_{l(k)+1}))-b_{s,
k}|-|w_k-a_{s, k}|
\end{equation}
$$\geqslant\delta/(m-1)-|f(\gamma_k(t_{l(k)+1})-b_{s,
k})|-|w_k-a_{s, k}|\,.$$
Since $a_{s, k}:=f(\gamma_k(p_s))\rightarrow
w_k:=f(\gamma_k(t_{l(k)}))$ as $s\rightarrow \infty$ and $b_{s,
k}=f(\gamma_k(q_s))\rightarrow f(\gamma_k(t_{l(k)+1}))$ as
$s\rightarrow \infty,$ it follows from the latter inequality that
there is $s=s(k)\in {\Bbb N}$ such that
\begin{equation}\label{eq12}
|a_{s(k), k}- b_{s(k), k}|\geqslant\delta/2(m-1)\,.
\end{equation}
Since $V$ is open, there is some neighborhood $U_k$ of $w_k$ such
$U_k\subset V.$ Since $a_{s, k}\rightarrow
w_k:=f(\gamma_k(t_{l(k)}))$ as $s\rightarrow \infty,$ we may assume
that $a_{s(k), k}\in V.$
Now, we set
$$\widetilde{G}_k:=f(\gamma_k)|_{[p_{s(k)}, q_{s(k)}]}\,.$$
In other words, $\widetilde{G}_k$ is a part of the path
$f(\gamma_k)$ between points $a_{s(k), k}$ and $b_{s(k), k}.$ Let us
to show that $\widetilde{G}_k$ intersects $\partial U$ and $\partial
V.$ Indeed, by the mentioned above, $a_{s(k), k}\subset V,$ so that
$V\cap \widetilde{G}_k\ne\varnothing.$ In particular, $U\cap
\widetilde{G}_k\ne\varnothing,$ because $V\subset U.$ On the other
hand, by (\ref{eq12}) $d(\widetilde{G}_k)\geqslant\delta/2(m-1),$
however, $d(U)< \delta/2(m-1)$ by the choice of $U.$ In particular,
$d(V)<\delta/2(m-1).$ It follows from this that
$$(\overline{{\Bbb R}^n}\setminus U)\cap
\widetilde{G}_k\ne\varnothing\,, \qquad(\overline{{\Bbb
R}^n}\setminus V)\cap \widetilde{G}_k\ne\varnothing\,.$$
Now, by \cite[Theorem~1.I.5, \S46]{Ku}
$$\partial U\cap
\widetilde{G}_k\ne\varnothing\,, \partial V\cap
\widetilde{G}_k\ne\varnothing\,,$$
as required.
Now, by~(\ref{eq1}),
\begin{equation}\label{eq1C}
M_{\alpha}(\Gamma(\widetilde{G}_k, C_0^{\,\prime}, K))\geqslant P
>0\,, \qquad k=1,2,\ldots\,.
\end{equation}
Let us to show that, the relation~(\ref{eq1C}) contradicts the
definition of $f$ and conditions~(\ref{eq7***})--(\ref{eq5***}).
Indeed, let us denote by $\Gamma_k$ the family of all half-open
paths $\beta_k:[a, b)\rightarrow \overline{{\Bbb R}^n}$ such that
$\beta_k(a)\in |f(\gamma_k)|,$ $\beta_k(t)\in K$ for all $t\in [a,
b)$ and, moreover, $\lim\limits_{t\rightarrow b-0}\beta_k(t):=B_k\in
C_0^{\,\prime}.$ Obviously, by~(\ref{eq1C})
\begin{equation}\label{eq4A}
M_{\alpha}(\Gamma_k)=M_{\alpha}(\Gamma(\widetilde{G}_k, F,
K))\geqslant P
>0\,, \qquad k=1,2,\ldots\,.
\end{equation}
Since by the proving above $D\setminus E$ consists of finite number
of components, we may assume that all paths
$\nabla_k:=\gamma_k|_{[p_{s(k)}, q_{s(k)}]}$ belong to (some) one
component $D_1$ of $D\setminus E.$ By the proving above, $f$ is
closed in $D_1$ and $f(D_1)=K.$ Consider the family
$\Gamma_k^{\,\prime}$ of maximal $f$-liftings $\alpha_k:[a,
c)\rightarrow D$ of the family $\Gamma_k$ starting at $|\nabla_k|;$
such a family exists by Proposition~\ref{pr3}.

Note that, $C(f, \partial D_1)\subset C(f, \partial D) \cup C(f,
E)\subset E_*.$ Now, observe that, the situation when
$\alpha_k\rightarrow
\partial D_1$ as $k\rightarrow\infty$ is impossible because $f$ is
closed in $D_1$ and, consequently, $f$ is boundary preserving (see
\cite[Theorem~3.3]{Vu}). Therefore, by Proposition~\ref{pr3}
$\alpha_k\rightarrow x_1\in D$ as $t\rightarrow c-0,$ and $c=b$ and
$f(\alpha_k(b))=f(x_1).$ In other words, the $f$-lifting $\alpha_k$
is complete, i.e., $\alpha_k:[a, b]\rightarrow D_1.$ Besides that,
it follows from that $\alpha_k(b)\in f^{\,-1}(C^{\,\prime}_0)\cap
D_1.$ Since $f(D_1)=K,$ we obtain that
$f(f^{\,-1}(C^{\,\prime}_0)\cap D_1)=C^{\,\prime}_0.$ In addition,
$f^{\,-1}(C^{\,\prime}_0)\cap D_1$ is a compactum in $D_1,$ see e.g.
\cite[Theorem~3.3]{Vu}. Thus, there is $r^{\,\prime}_0>0$ such that
\begin{equation}\label{eq13}
d(f^{\,-1}(C^{\,\prime}_0)\cap D_1, \partial D)\geqslant r_0>0\,.
\end{equation}
Moreover, it follows from~(\ref{eq13}) that
\begin{equation}\label{eq14}
f^{\,-1}(C^{\,\prime}_0)\subset D\setminus B(b, r_0)\,.
\end{equation}
We may consider that $r_0<\varepsilon_0,$ where $\varepsilon_0$ is a
number from the conditions of the lemma. Applying now
Lemmas~\ref{thOS4.1} and~\ref{lem4}, we obtain for
$\alpha:=\frac{p}{p-n+1}$ that
\begin{gather}
M_{\alpha}(f(\Gamma_k^{\,\prime}))=M_{\alpha}(\Gamma(\widetilde{G}_k,
C_0^{\,\prime}, K)) \nonumber \\ \label{eq15}\leqslant(N(f,
D_1))^{\frac{n-1}{p-n+1}}\int\limits_{A(b, 2^{\,-k},
r_0)}Q(x)\eta^{\,\alpha}(|x-b|)\,dm(x)
\end{gather}
for sufficiently large $k\in {\Bbb N}$ for any nonnegative Lebesgue
measurable function $\eta$ satisfying the relation~(\ref{eq6B}) with
$\varepsilon=2^{\,-k}$ and $\varepsilon_0=r_0.$ It is worth noting
that $N(f, D_1)<\infty$ because $f$ is closed in $D_1,$ see e.g.
\cite[Lemma~3.3]{MS}. Observe that the function
$$\eta(t)=\left\{
\begin{array}{rr}
\psi(t)/I(2^{\,-k}, r_0), &   t\in (2^{\,-k},
r_0),\\
0,  &  t\in {\Bbb R}\setminus (2^{\,-k}, r_0)\,,
\end{array}
\right. $$
where $I(2^{\,-k}, r_0)$ is defined in~(\ref{eq7***}),
satisfies~(\ref{eq6B}). Therefore, by~(\ref{eq15}) and
(\ref{eq7***})--(\ref{eq5***}) we obtain that
\begin{equation}\label{eq11*}
M_{\alpha}(f(\Gamma^{\,\prime}_k))=M_{\alpha}(\Gamma(\widetilde{G}_k,
C_0^{\,\prime}, K))\leqslant (N(f,
D_1))^{\frac{n-1}{p-n+1}}\Delta(k)\,,
\end{equation}
where $\Delta(k)\rightarrow 0$ as $k\rightarrow \infty.$ The
relations~(\ref{eq11*}) together with~(\ref{eq1C}) contradict each
other, which proves the lemma. $\Box$
\end{proof}

\medskip
{\it Proof of Theorem~\ref{th3}} immediately follows by
Lemma~\ref{lem1} and Proposition~\ref{pr6}, excluding the equality
$f(\overline{D})=\overline{D^{\,\prime}}.$ The proof of the latter
repeats the last part of the proof of Theorem~3.1 in
\cite{SSD}.~$\Box$

\medskip
{\it Proof of Theorem~\ref{th1}} is similar to the proof of
Theorem~\ref{th3}, but uses Lemma~\ref{lem4} instead of
Lemma~\ref{thOS4.1}. In all other respects, the proof of this
theorem literally repeats the previous proof.~$\Box$

A domain $D\subset {\Bbb R}^n,$ $n\geqslant 2,$ is called a {\it
uniform} domain with respect to $p$-modulus, $p\geqslant 1,$ if, for
each $r>0,$ there is $\delta>0$ such that $M_{p}(\Gamma(F, F^{\,*},
D))\geqslant\delta$ whenever $F$ and $F^{\,*}$ are continua of $D$
with $h(F)\geqslant r$ and $h(F^{\,*})\geqslant r.$ Let $I$ be some
set of indices. Domains $D_i,$ $i\in I,$ are said to be {\it
equi-uniform} domains with respect to $p$-modulus, if, for $r>0,$
the modulus condition above is satisfied by each $D_i$ with the same
number $\delta.$ It should be noted that the proposed concept of a
uniform domain has, generally speaking, no relation to definition,
introduced for the uniform domain in Martio-Sarvas sense~\cite{MSa}.
Note that, uniform domains with respect to $p$-modulus have strongly
accessible boundaries with respect to $p$-modulus, see
\cite[Remark~1]{SevSkv$_1$}.

\medskip
\begin{remark}\label{rem2}
The statement of Lemma~\ref{lem1} remains true if condition 3) in
this lemma is replaced by the conditions:

($3_\textbf{a}$) the family of components of $D^{\,\prime}\setminus
E_*$ is {\it equi-uniform} with respect to $p$-modulus and, besides
that,

(\textbf{$3_\textbf{b}$}) there is $\delta_*>0$ such that for each
component $K$ of $D^{\,\prime}\setminus  E_*$ there is a
nondegenerate continuum $F\subset K$ such that $d(F)\geqslant
\delta_*$ and $d(f^{\,-1}(F),
\partial D)\geqslant \delta_*>0.$
\end{remark}

\medskip
The proof of Remark~\ref{rem2} may be given similarly to Remark~4.1
in~\cite{DS}.

\section{Some additional information}

Later, in the extended space $\overline{{{\Bbb R}}^n}={{\Bbb
R}}^n\cup\{\infty\}$ we use the {\it spherical (chordal) metric}
$h(x,y)=|\pi(x)-\pi(y)|,$ where $\pi$ is a stereographic projection
$\overline{{{\Bbb R}}^n}$ onto the sphere
$S^n(\frac{1}{2}e_{n+1},\frac{1}{2})$ in ${{\Bbb R}}^{n+1},$ namely,
\begin{equation}\label{eq3C}
h(x,\infty)=\frac{1}{\sqrt{1+{|x|}^2}}\,,\quad
h(x,y)=\frac{|x-y|}{\sqrt{1+{|x|}^2} \sqrt{1+{|y|}^2}}\,, \quad x\ne
\infty\ne y
\end{equation}
(see \cite[Definition~12.1]{Va}). Further, the closure
$\overline{A}$ and the boundary $\partial A$ of the set $A\subset
\overline{{\Bbb R}^n}$ we understand relative to the chordal metric
$h$ in $\overline{{\Bbb R}^n}.$

\medskip
Following \cite[Section~II.10]{Ri}, a {\it condenser} is a pair
$E=(A, C)$ where $A\subset {\Bbb R}^n$ is open and $C$ is non--empty
compact set contained in $A.$ Given $p\geqslant 1$ and a condenser
$E=\left(A,\,C\right),$ we set
%
%\begin{equation}\label{equ5}
$${\rm cap}_p\,E\quad=\quad{\rm
cap}_p\,\left(A,\,C\right)\quad=\quad\inf \limits_{u\in
W_0\left(E\right) }\quad\int\limits_{A}\,|\nabla u|^p dm(x)$$
%\end{equation}
%
%
where $W_0(E)=W_0(A,\,C)$ is a family of all non-negative functions
$u:A\rightarrow {\Bbb R}^1$ such that (1)\quad $u\in
C_0(A),$\quad(2)\quad $u(x)\ge 1$ for $x\in C,$ and (3)\quad $u$ is
$ACL.$ In the above formula
$|\nabla u|={\left(\sum\limits_{i=1}^n\,{\left(\partial_i
u\right)}^2 \right)}^{1/2},$ and ${\rm cap}_p\,E$ is called {\it
$p$-ca\-pa\-ci\-ty} of the condenser $E,$ see in \cite[Section
~II.10]{Ri}.

\medskip
The following lemma was proved in \cite[Lemma~4.4]{Sev$_1$},
cf.~\cite[Lemma~4.2]{SevSkv$_3$}.

\medskip
\begin{lemma}\label{lem1A} {\it Let $x_0\in D,$ let $n\geqslant 2,$ let $p>n-1$ and let
$Q:D\rightarrow [0,\infty]$ be locally integrable function with
degree $\alpha:=\frac{n-1}{p-n+1}.$ If $f:D\rightarrow
\overline{{\Bbb R}^n}$ is an open discrete lower $Q$-mapping at
$x_0$ with respect to $p$-modulus, then $f$ satisfies
$${\rm cap}_{\beta}(f(E))\leqslant \frac{\omega_{n-1}}{I^{\,*\,\beta-1}}$$
at the point $x_0$ with $\beta:=\frac{p}{p-n+1}$ and
$Q^{\,*}=Q^{\frac{n-1}{p-n+1}}.$

Here $E$ is the condenser $(B(x_0, r_2), \overline{B(x_0, r_1)}),$
$q^{\,*}_{x_0}(r)$ denotes the spherical mean of
$Q^{\frac{n-1}{p-n+1}}$ over $|x-x_0|=r$ and $I^{\,*}=I^{\,*}(r_1,
r_2)=\int\limits_{r_1}^{r_2}\frac{dt}{t^{\frac{n-1}{\beta-1}}
q_{x_0}^{*\,\frac{1}{\beta-1}}(t)}$ for $0<r_1<r_2<{\rm dist\,}(x_0,
\partial D).$
}\end{lemma}

\medskip In what follows we will need the following auxiliary
assertion (see, for example, \cite[Lemma~7.4, ch.~7]{MRSY$_2$} for
$p=n$ and \cite[Lemma~2.2]{Sal} for $p\ne n.$

\medskip
\begin{proposition}\label{pr1A}
{\,\it Let $x_0 \in {\Bbb R}^n,$ $Q(x)$ be a Lebesgue measurable
function, $Q:{\Bbb R}^n\rightarrow [0, \infty],$ $Q\in
L_{loc}^1({\Bbb R}^n).$ We set $A:=A(x_0, r_1, r_2)=\{ x\,\in\,{\Bbb
R}^n : r_1<|x-x_0|<r_2\}$ and
$\eta_0(r)=\frac{1}{Ir^{\frac{n-1}{p-1}}q_{x_0}^{\frac{1}{p-1}}(r)},$
where $I:=I=I(x_0,r_1,r_2)=\int\limits_{r_1}^{r_2}\
\frac{dr}{r^{\frac{n-1}{p-1}}q_{x_0}^{\frac{1}{p-1}}(r)}$ and
$q_{x_0}(r):=\frac{1}{\omega_{n-1}r^{n-1}}\int\limits_{|x-x_0|=r}Q(x)\,d{\mathcal
H}^{n-1}$ is the integral average of the function $Q$ over the
sphere $S(x_0, r).$ Then
\begin{equation*}\label{eq10A_1}
\frac{\omega_{n-1}}{I^{p-1}}=\int\limits_{A} Q(x)\cdot
\eta_0^p(|x-x_0|)\ dm(x)\leqslant\int\limits_{A} Q(x)\cdot
\eta^p(|x-x_0|)\ dm(x)
\end{equation*}
for any Lebesgue measurable function $\eta :(r_1,r_2)\rightarrow
[0,\infty]$ such that
$\int\limits_{r_1}^{r_2}\eta(r)\,dr=1. $ Moreover, (\ref{eq10A_1})
holds for similar functions $\eta$ with
$\int\limits_{r_1}^{r_2}\eta(r)\,dr\geqslant 1$ (see, e.g.,
\cite[Remark~3.1]{Sev$_4$}).}
\end{proposition}

\medskip
Combining Lemma~\ref{lem1A} with Proposition~\ref{pr1A}, we obtain
the following statement.

\medskip
\begin{lemma}\label{lem5} {\it Let $x_0\in D,$ let $n\geqslant 2,$ let $p>n-1$ and let
$Q:D\rightarrow [0,\infty]$ be locally integrable function with
degree $\alpha:=\frac{n-1}{p-n+1}.$ If $f:D\rightarrow
\overline{{\Bbb R}^n}$ is an open discrete lower $Q$-mapping at
$x_0$ with respect to $p$-modulus, $p>n-1,$ then $f$ satisfies
$${\rm cap}_{\beta}(f(E))\leqslant \int\limits_{A(x_0, r_1, r_2)} Q^{\frac{n-1}{p-n+1}}(x)\cdot
\eta^{\beta}(|x-x_0|)\ dm(x)$$
for any Lebesgue measurable function $\eta :(r_1,r_2)\rightarrow
[0,\infty]$ such that
$\int\limits_{r_1}^{r_2}\eta(r)\,dr\geqslant 1,$ where
$\beta:=\frac{p}{p-n+1}$ and $Q^{\,*}=Q^{\frac{n-1}{p-n+1}}.$

Here $E$ is the condenser $(B(x_0, r_2), \overline{B(x_0, r_1)}),$
$0<r_1<r_2<{\rm dist\,}(x_0,
\partial D).$
}\end{lemma}

\medskip
The following definition is from \cite[Section~2, Ch.~III]{Ri}. Let
$F$ be a compact set in ${\Bbb R}^n$. We say that $F$ is of {\it
$p$-capacity zero} if ${\rm cap}_p\,(A, F)=0$ for some bounded open
set $A\supset F.$ An arbitrary set $E\subset {\Bbb R}^n$ is of
$p$-capacity zero if the same is true for every compact subset of
$E.$ In this case we write ${\rm cap}_p\, E = 0$ (${\rm cap}\, E=0$
if $p=n$), otherwise ${\rm cap}_p E>0.$ The following statement was
proved in~\cite[Lemma~4]{Sev$_3$} for $p=n$.

\medskip
\begin{lemma}\label{lem6}
{\it\, Let $n-1<p\leqslant n,$ let $D$ be a domain in ${\Bbb R}^n,$
$n\ge 2,$ let $E\subset \overline{{\Bbb R}^n}$ be a compact set of
positive capacity for $p=n,$ and an arbitrary set for $p\ne n,$ let
$Q:D\rightarrow [0, \infty]$ be a Lebesgue measurable function, and
let ${\frak F}^{p}_{Q, E}(x_0)$ be a family of open discrete
mappings $f:D\,\rightarrow\,\overline{{\Bbb R}^n}\setminus E$
satisfying condition
$${\rm cap}_{p}(f({\mathcal E}))\leqslant \int\limits_{A(x_0, r_1, r_2)} Q(x)\cdot
\eta^p(|x-x_0|)\ dm(x)$$
for any Lebesgue measurable function $\eta :(r_1,r_2)\rightarrow
[0,\infty]$ such that
$\int\limits_{r_1}^{r_2}\eta(r)\,dr\geqslant 1,$ ${\mathcal
E}=(B(x_0, r_2), \overline{B(x_0, r_1)}).$
Assume that, the relations~(\ref{eq7***})--(\ref{eq5***}) holds with
$b=x_0.$ Then the family of mappings ${\frak F}^{p}_{Q, E}(x_0)$ is
equicontinuous at the point~$x_0.$ }
\end{lemma}

\medskip
\begin{proof}
As we have already noted, this statement was proved earlier for
$p=n,$ see \cite[Lemma~4]{Sev$_3$}. The proof for the case $p\ne n$
literally repeats the proof of Lemma~4.3 in \cite{Sev$_1$} and is
therefore omitted.~$\Box$
\end{proof}

\section{The equicontinuity in the closure of a domain}

Let us formulate some statements about the equicontinuity of the
Sobolev and Orlicz-Sobolev classes in the closure of the domain. For
some of our other results on this topic, see, for example,
\cite{SevSkv$_1$}, \cite{SevSkv$_2$} and \cite{DS}.

\medskip
Given $N\in {\Bbb N},$ $\delta>0,$ $\alpha\geqslant 1,$ closed sets
$E_*, F$ in $\overline{{\Bbb R}^n},$ $n\geqslant 2,$ a domain
$D\subset {\Bbb R}^n,$ a closed (with respect to $D$) set $E$ in
$D,$ an increasing function
$\varphi:[0,\infty)\rightarrow[0,\infty)$ and a Lebesgue measurable
function $Q:D\rightarrow [0, \infty]$ let us denote by
$\frak{R}^{E_*, E, F, N}_{\varphi, Q, \delta, \alpha}(D)$ a (some)
family of bounded open discrete mappings $f:D\rightarrow
\overline{{\Bbb R}^n}\setminus F$ such that $K_{I, \alpha}(x,
f)\leqslant Q(x)$ a.e., $N(f, D)\leqslant N$ and, in addition,

\medskip
1) $C(f, \partial D)\subset E_*,$

\medskip
2) for each component $K$ of $D^{\,\prime}_f\setminus  E_*,$
$D^{\,\prime}_f:=f(D),$ there is a continuum $K_f\subset K$ such
that $h(K_f)\geqslant \delta$ and $h(f^{\,-1}(K_f), \partial
D)\geqslant \delta>0,$

\medskip
3) $f^{\,-1}(E_*)=E.$

\medskip
The following statement holds.

\medskip
\begin{theorem}\label{th4} {\it\, Let $\alpha\in(n-1, n],$
let $D$ be a bounded domain in ${\Bbb R}^n,$ $n\geqslant 2.$
Assume that:

\medskip
1) the set $E$ is nowhere dense in $D,$ and $D$ is finitely
connected on $E\cup \partial D,$ i.e., for any $z_0\in E\cup
\partial D$ and any neighborhood $\widetilde{U}$ of $z_0$ there is a
neighborhood $\widetilde{V}\subset \widetilde{U}$ of $z_0$ such that
$(D\cap \widetilde{V})\setminus E$ consists of finite number of
components;

\medskip
2) for any $x_0\in\partial D$ there is $m=m(x_0)\in {\Bbb N},$
$1\leqslant m<\infty$ such that the following is true: for any
neighborhood $U$ of $x_0$ there is a neighborhood $V\subset U$ of
$x_0$ and  such that:

\medskip
2a) $V\cap D$ is connected,

\medskip
2b) $(V\cap D)\setminus E$ consists at most of $m$ components.

\medskip
Let for $\alpha=n$ the set $F$ have positive capacity, and for
$n-1<\alpha<n$ it is an arbitrary closed set.

\medskip
Suppose that, for any $x_0\in\overline{D}$ at least one of the
following conditions is satisfied: $3_1)$ a function $Q$ has a
finite mean oscillation at $x_0;$ $3_2)$
$q_{x_0}(r)\,=\,O\left(\left[\log{\frac1r}\right]^{n-1}\right)$ as
$r\rightarrow 0;$ $3_3)$ the condition
\begin{equation*}\label{eq6D}
\int\limits_{0}^{\delta(x_0)}\frac{dt}{t^{\frac{n-1}{\alpha-1}}
q_{x_0}^{\,\prime\,\frac{1}{\alpha-1}}(t)}=\infty
\end{equation*}
holds for some $\delta(x_0)>0,$ where $q_{x_0}^{\,\prime}(t)$ is
defined in~(\ref{eq32*B}).

\medskip
4) Assume that, the function $\varphi$ satisfies Calderon
condition~(\ref{eq1_A_10}).

Let the family of all components of $D^{\,\prime}_f\setminus E_*$ is
equi-uniform over $f\in\frak{R}^{E_*, E, F, N}_{\varphi, Q, \delta,
\alpha}(D)$ with respect to $\alpha$-modulus. Then every
$f\in\frak{R}^{E_*, E, F, N}_{\varphi, Q, \delta, \alpha}(D)$ has a
continuous extension to $\partial D$ and the family $\frak{R}^{E_*,
E, F, N}_{\varphi, Q, \delta, \alpha}(\overline{D}),$ consisting of
all extended mappings $\overline{f}: \overline{D}\rightarrow
\overline{{\Bbb R}^n},$ is equicontinuous in~$\overline{D}.$ The
equicontinuity must be understood in the sense of the chordal metric
$h$ defined in~(\ref{eq3C}).}
\end{theorem}

\medskip
The proof of Theorem~\ref{th3} is based on the following lemma.

\medskip
\begin{lemma}\label{lem3} {\it\, The conclusion of Theorem~\ref{th4}
remains valid if, under conditions of this theorem, we replace the
assumptions $3_1)$--$3_3)$ by the following assumption: Suppose
that, for any $x_0\in \overline{D}$ there is
$\varepsilon_0=\varepsilon_0(x_0)>0$ and some positive measurable
function $\psi:(0, \varepsilon_0)\rightarrow (0,\infty)$ such that
\begin{equation}\label{eq8}
0<I(\varepsilon,
\varepsilon_0)=\int\limits_{\varepsilon}^{\varepsilon_0}\psi(t)\,dt
< \infty
\end{equation}
for sufficiently small $\varepsilon\in(0, \varepsilon_0)$ and, in
addition,
\begin{equation}\label{eq9}
\int\limits_{A(x_0, \varepsilon, \varepsilon_0)}
Q(x)\cdot\psi^{\,\alpha}(|x-x_0|)
 \ dm(x) =o(I^{\,\alpha}(\varepsilon, \varepsilon_0))\,,
\end{equation}
where $A:=A(x_0, \varepsilon, \varepsilon_0)$ is defined
in~(\ref{eq1**}). }
\end{lemma}

\medskip
\begin{proof}
Let $x_0\in D.$ By Lemma~\ref{thOS4.1}, every $f\in \frak{R}^{E_*,
E, F, N}_{\varphi, Q, \delta, \alpha}(D)$ is a lower $Q$-mapping
relative to the $p$-modulus at every point $x_0\in\overline {D}$ for
$Q(x)=N\cdot K^{\frac{p-n+1}{n-1}}_{I, \alpha}(x, f),$ where
$\alpha=\frac{p}{p-n+1}.$ Now, by Lemma~\ref{lem5} $f$ satisfies the
relation
$${\rm cap}_{\alpha}(f({\mathcal{E}}))\leqslant \int\limits_{A(x_0,
r_1, r_2)} N^{\frac{n-1}{p-n+1}}\cdot K_{I, \alpha}(x,
f)\eta^{\,\alpha}(|x-x_0|)\,dm(x)$$
$$\leqslant
 N^{\frac{n-1}{p-n+1}}\cdot\int\limits_{A(x_0, r_1, r_2)} Q(x)\eta^{\,\alpha}(|x-x_0|)\,dm(x)$$
for any Lebesgue measurable function $\eta :(r_1,r_2)\rightarrow
[0,\infty]$ such that
$\int\limits_{r_1}^{r_2}\eta(r)\,dr\geqslant 1.$ Now, the
equicontinuity of the family $\frak{R}^{E_*, E, F, N}_{\varphi, Q,
\delta, \alpha}(D)$ at inner points $x_0\in D$ follows by
Lemma~\ref{lem6}.

\medskip
Since any component of $D^{\,\prime}_f\setminus E_*$ is uniform with
respect to $\alpha$-modulus, it has strongly accessible boundary
with respect to $\alpha$-modulus (see \cite[Remark~1]{SevSkv$_1$}).
Now the continuous extension of any $f\in\frak{R}^{E_*, E, F,
N}_{\varphi, Q, \delta, \alpha}(D)$ to $\partial D$ follows from
Lemma~\ref{lem1} and Remark~\ref{rem2}. It remains to prove the
equicontinuity of the extended family $\frak{R}^{E_*, E, F,
N}_{\varphi, Q, \delta, \alpha}(\overline{D})$ in~$\partial D.$

Suppose the opposite. Then there is $x_0\in\partial D,$
$\varepsilon_0>0,$ a sequence $x_k\rightarrow x_0,$
$x_k\in\overline{D},$ and $f_k\in \frak{R}^{E_*, E, F, N}_{\varphi,
Q, \delta, \alpha}(D)$ such that
\begin{equation}\label{eq12A}
h(f_k(x_k), f_k(x_0))\geqslant \varepsilon_0\,.
\end{equation}
Since $f_k$ has a continuous extension to $\partial D,$ we may
consider that $x_k\in D.$ In addition, it follows from~(\ref{eq12A})
that we may find a sequence $x^{\,\prime}_k\in D,$ $k=1,2,\ldots ,$
such that $x^{\,\prime}_k\rightarrow x_0$ and
\begin{equation}\label{eq13A}
h(f_k(x_k), f_k(x^{\,\prime}_k))\geqslant \varepsilon_0/2\,.
\end{equation}
By the assumption~2), there exists a sequence of neighborhoods
$V_k\subset B(x_0, 2^{\,-k}),$ $k=1,2,\ldots ,$ such that $V_k\cap
D$ is connected and $(V_k\cap D)\setminus E$ consists of $m$
components, $1\leqslant m<\infty.$

Arguing similarly to the proof of Lemma~\ref{lem1}, we may consider
that $x_k, x^{\,\prime}_k\not\in E.$ Now, by Lemma~\ref{lem2} there
are subsequences $x_{k_l}$ and $x^{\,\prime}_{k_l},$ $l=1,2,\ldots
,$ belonging to some sequence of neighborhoods $V_l,$ $l=1,2,\ldots
,$ of the point $x_0$ such that ${\rm diam\,}V_l\rightarrow 0$ as
$l\rightarrow\infty$ and, in addition, any pair $x_{k_l}$ and
$x^{\,\prime}_{k_l}$ may be joined by a path $\gamma_l$ in $V_l\cap
D,$ where $\gamma_l$ contains at most $m-1$ points in $E.$ Without
loss of generality, we may assume that the same sequences $x_k$ and
$y_k$ satisfy properties mentioned above. Let $\gamma_k:[0,
1]\rightarrow D,$ $\gamma_k(0)=x_k$ and
$\gamma_k(1)=x^{\,\prime}_k,$ $k=1,2,\ldots .$

\medskip
Observe that, the path $f_k(\gamma_k)$ contains not more than $m-1$
points in $E_*.$ Let
$$b_{1}=f_k(\gamma_k(t_1)), b_{2}=f_k(\gamma_k(t_2))\quad,\ldots,\quad
b_l=f_k(\gamma_k(t_l))\,,$$ $$t_0:=0\leqslant t_1\leqslant
t_2\leqslant \ldots\leqslant t_l\leqslant 1=t_{l+1},\qquad
0\leqslant l\leqslant m-1\,,$$
be points in $f_k(\gamma_k)\cap E_*.$ By the relation~(\ref{eq13A})
and due to the triangle inequality,
\begin{equation}\label{eq6C}
\varepsilon_0/2\leqslant h(f_k(x_k),
f_k(x^{\,\prime}_k))\leqslant\sum\limits_{r=0}^{l}
h(f_k(\gamma_k(t_r)), f_k(\gamma_k(t_{r+1}))\,.
\end{equation}
It follows from~(\ref{eq6C}) that, there is $1\leqslant
r=r(k)\leqslant m-1$ such that
\begin{equation}\label{eq7F}
h(f(\gamma_k(t_{r(k)})), f_k(\gamma_k(t_{r(k)+1})))\geqslant
\varepsilon_0/(2(l+1))\geqslant \varepsilon_0/(2m)\,.
\end{equation}
Observe that, the set $G_k:=|f(\gamma_k)|_{(t_{r(k)}, t_{r(k)+1})}|$
belongs to $D^{\,\prime}\setminus E_*$ and that $G_k$ contains some
a continuum $\widetilde{G}_k$ with $h(\widetilde{G}_k)\geqslant
\varepsilon_0/(4m)$ for any $k\in {\Bbb N}.$

Since $\widetilde{G}_k$ is a continuum in
$D^{\,\prime}_{f_k}\setminus E_*,$ there is a component $K_k$ of
$D^{\,\prime}_{f_k}\setminus E_*,$ containing $K_k.$ Let us apply
the definition of equi-uniformity for he sets $\widetilde{G}_k$ and
$K_{f_k}$ in $K_k$ (here $K_{f_k}$ is a continuum from the
definition of the class $\frak{R}^{E_*, E, F, N}_{\varphi, Q,
\delta, \alpha}(D)$, in particular, $h(K_{f_k})\geqslant \delta$).
Due to this definition, for the number $\delta_*:=\min\{\delta,
\varepsilon/4m\}>0$ there is $P>0$ such that
\begin{equation}\label{eq1F}
M_p(\Gamma(\widetilde{G}_k, K_{f_k}, K_k))\geqslant P
>0\,, \qquad k=1,2,\ldots\,.
\end{equation}
Let us to show that, the relation~(\ref{eq1F}) contradicts with the
definition of the class $\frak{R}^{E_*, E, F, N}_{\varphi, Q,
\delta, \alpha}(D)$ together with
conditions~(\ref{eq7***})--(\ref{eq5***}). Indeed, let us denote by
$\Gamma_k$ the family of all half-open paths $\beta_k:[a,
b)\rightarrow \overline{{\Bbb R}^n}$ such that $\beta_k(a)\in
\widetilde{G}_k,$ $\beta_k(t)\in K_k$ for all $t\in [a, b)$ and,
moreover, $\lim\limits_{t\rightarrow b-0}\beta_k(t):=B_k\in
K_{f_k}.$ Obviously, by~(\ref{eq1F})
\begin{equation}\label{eq4C}
M_p(\Gamma_k)=M_p(\Gamma(\widetilde{G}_k, K_{f_k}, K_k))\geqslant P
>0\,, \qquad k=1,2,\ldots\,.
\end{equation}
As under the proof of Lemma~\ref{lem1}, we may prove that
$D\setminus E$ consists of finite number of components $D_1,\ldots,
D_s,$ $1\leqslant s<\infty,$ while $f_k$ is closed in each $D_i,$
$i=1,2,\ldots, $ i.e. $f_k(E)$ is closed in $f_k(D_i)$ whenever $E$
is closed in $D_i.$ We also may show that $f_k(D_i)=K$ for some
$i\in {\Bbb N}.$ Without loss of generality, we may consider that
$\nabla_k:=\gamma_k|_{[p_{s(k)}, q_{s(k)}]}$ belong to (some) one
component $D_1$ of $D\setminus E$ and $f_k(D_1)=K,$ where
$t_{r(k)}<p_{s(k)}<q_{s(k)}<t_{r(k)+1},$ $a_{s,
k}:=f_k(\gamma_k(p_s))\rightarrow w_k:=f_k(\gamma_k(t_{r(k)}))$ as
$s\rightarrow \infty$ and $b_{s, k}=f_k(\gamma_k(q_s))\rightarrow
f_k(\gamma_k(t_{r(k)+1}))$ as $s\rightarrow \infty.$ Consider the
family $\Gamma_k^{\,\prime}$ of all maximal $f_k$-liftings
$\alpha_k:[a, c)\rightarrow D$ of the family $\Gamma_k$ starting at
$|\nabla_k|;$ such a family exists by Proposition~\ref{pr3}.

Observe that, the situation when $\alpha_k\rightarrow \partial D_1$
as $k\rightarrow\infty$ is impossible. Suppose the opposite: let
$\alpha_k(t)\rightarrow \partial D_1$ as $t\rightarrow c.$ We choose
an arbitrary sequence $\varphi_m\in [0, c)$ such that
$\varphi_m\rightarrow c-0$ as $m\rightarrow\infty.$ Since the space
$\overline{{\Bbb R}^n}$ is compact, the boundary $\partial D_1$ is
also compact as a closed subset of the compact space. Then there
exists $w_m\in
\partial D_1$ such that
\begin{equation}\label{eq7G}
h(\alpha_k(\varphi_m), \partial D_1)=h(\alpha_k(\varphi_m), w_m)
\rightarrow 0\,,\qquad m\rightarrow \infty\,.
\end{equation}
Due to the compactness of $\partial D_1$, we may assume that
$w_m\rightarrow w_0\in \partial D_1$ as $m\rightarrow\infty.$
Therefore, by the relation~(\ref{eq7G}) and by the triangle
inequality
\begin{equation}\label{eq8C}
h(\alpha_k(\varphi_m), w_0)\leqslant h(\alpha_k(\varphi_m),
w_m)+h(w_m, w_0)\rightarrow 0\,,\qquad m\rightarrow \infty\,.
\end{equation}
There are two cases: $w_0\in \partial D,$ or $w_0\in D.$ In the
first case,
\begin{equation}\label{eq9C}
f_k(\alpha_k(\varphi_m))=\beta_k(\varphi_m)\rightarrow \beta(c)
\,,\quad m\rightarrow\infty\,,
\end{equation}
because by the construction the path $\beta_k(t),$ $t\in [a, b],$
lies in $K_k\subset D^{\,\prime}_{f_k}\setminus E_*$ together with
its finite ones points. At the same time, by~(\ref{eq8C})
and~(\ref{eq9C}) we have that $\beta_k(c)\in C(f_k,
\partial D)\subset E_*$ by the definition of the class
$\frak{R}^{E_*, E, F, N}_{\varphi, Q, \delta, \alpha}(D)$ and
because $C(f_k, \partial D)\subset E_*.$ The inclusions
$\beta_k\subset D^{\,\prime}\setminus E_*$ and $\beta_k(c)\in E_*$
contradict each other.

\medskip
In the second case, when $w_0\in D,$ we have that $w_0\in E.$
However,
$f_k(\alpha_k(\varphi_m))=\beta_k(\varphi_m)\rightarrow\beta_k(c)$
and consequently, $\beta_k(c)\in f_k(E).$ But now $\beta_k(c)\in
E_*$ because $f_k^{\,-1}(E_*)=E.$ The latter contradicts the
definition of $\beta_k.$ Thus, the case $\alpha_k\rightarrow
\partial D_1$ as $k\rightarrow\infty$ is impossible, as required.

\medskip
Therefore, by Proposition~\ref{pr3} $\alpha_k\rightarrow x_1\in D_1$
as $t\rightarrow c-0,$ and $c_b$ and $f_k(\alpha_k(b))=f_k(x_1).$ In
other words, the $f_k$-lifting $\alpha_k$ is complete, i.e.,
$\alpha_k:[a, b]\rightarrow D.$ Besides that, it follows from that
$\alpha_k(b)\in f_k^{\,-1}(K_{f_k}).$

\medskip
Again, by the definition of the class $\frak{R}^{E_*, E, F,
N}_{\varphi, Q, \delta, \alpha}(D),$
\begin{equation}\label{eq13B}
h(f_k^{\,-1}(K_{f_k}), \partial D)\geqslant \delta>0\,.
\end{equation}
Since $x_0\ne\infty,$ it follows from~(\ref{eq13B}) that
\begin{equation}\label{eq14C}
f_k^{\,-1}(K_{f_k})\subset D\setminus B(x_0, r_0)
\end{equation}
for any $k\in {\Bbb N}$ and some $r_0>0.$ Let $k$ be such that
$2^{\,-k}<\varepsilon_0.$ We may consider that $r_0<\varepsilon_0,$
where $\varepsilon_0$ is a number from the conditions of the lemma.
Applying now Lemmas~\ref{thOS4.1} and~\ref{lem4}, we obtain for
$\alpha:=\frac{p}{p-n+1}$ that
\begin{equation}\label{eq15A}
M_{\alpha}(f(\Gamma_k^{\,\prime}))=M_{\alpha}(\Gamma(\widetilde{G}_k,
C_0^{\,\prime}, K))\leqslant
N^{\frac{n-1}{p-n+1}}\int\limits_{A(x_0, 2^{\,-k},
r_0)}Q(x)\eta^{\,\alpha}(|x-x_0|)\,dm(x)
\end{equation}
for sufficiently large $k\in {\Bbb N}$ and for any nonnegative
Lebesgue measurable function $\eta$ satisfying the
relation~(\ref{eq6B}) with $\varepsilon=2^{\,-k}$ and
$\varepsilon_0=r_0.$ Observe that the function
$$\eta(t)=\left\{
\begin{array}{rr}
\psi(t)/I(2^{\,-k}, r_0), &   t\in (2^{\,-k},
r_0),\\
0,  &  t\in {\Bbb R}\setminus (2^{\,-k}, r_0)\,,
\end{array}
\right. $$
where $I(2^{\,-k}, r_0)$ is defined in~(\ref{eq7***}),
satisfies~(\ref{eq6B}) for $\varepsilon:=2^{\,-k},$
$\varepsilon_0:=r_0.$ Therefore, by~(\ref{eq15A}) and
(\ref{eq7***})--(\ref{eq5***}) we obtain that
\begin{equation}\label{eq11*A}
M_{\alpha}(f(\Gamma^{\,\prime}_k))=M_{\alpha}(\Gamma(\widetilde{G}_k,
C_0^{\,\prime}, K))\leqslant N^{\frac{n-1}{p-n+1}}\cdot\Delta(k)\,,
\end{equation}
where $\Delta(k)\rightarrow 0$ as $k\rightarrow \infty.$ The
relations~(\ref{eq11*A}) together with~(\ref{eq1F}) contradict each
other, which proves the lemma. $\Box$
\end{proof}

\medskip {\it Proof of Theorem~\ref{th4}} directly follows from
Lemmas~\ref{lem3} and Proposition~\ref{pr6}.~$\Box$

\medskip
We also may formulate the corresponding statement for Sobolev space
on the plane, cf.~Theorem~\ref{th1}.

\medskip
Given $N\in {\Bbb N},$ $\delta>0,$ $\alpha\geqslant 1,$ closed sets
$E, E_*, F$ in $\overline{{\Bbb R}^n},$ $n\geqslant 2,$ a domain
$D\subset {\Bbb C}$ and a Lebesgue measurable function
$Q:D\rightarrow [0, \infty]$ let us denote by $\frak{R}^{E_*, E, F,
N}_{Q, \delta, \alpha}(D)$ the family of all bounded open discrete
mappings $f:D\rightarrow \overline{\Bbb C}\setminus F$ such that
$K_{I, \alpha}(x, f)\leqslant Q(x)$ a.e., $N(f, D)\leqslant N$ and,
in addition,

\medskip
1) $C(f, \partial D)\subset E_*,$

\medskip
2) for each component $K$ of $D^{\,\prime}_f\setminus  E_*,$
$D^{\,\prime}_f:=f(D),$ there is a continuum $K_f\subset K$ such
that $h(K_f)\geqslant \delta$ and $h(f^{\,-1}(K_f), \partial
D)\geqslant \delta>0,$

\medskip
3) $f^{\,-1}(E_*)=E.$

\medskip
The following statement holds.

\medskip
\begin{theorem}\label{th5} {\it\, Let $\alpha\in(1, 2],$
let $D$ be a bounded domain in ${\Bbb C}.$ Assume that:

\medskip
1) the set $E$ is nowhere dense in $D,$ and $D$ is finitely
connected on $E\cup\partial D,$ i.e., for any $z_0\in E\cup\partial
D$ and any neighborhood $\widetilde{U}$ of $z_0$ there is a
neighborhood $\widetilde{V}\subset \widetilde{U}$ of $z_0$ such that
$(D\cap \widetilde{V})\setminus E$ consists of finite number of
components;

\medskip
2) for any $x_0\in\partial D$ there is $m=m(x_0)\in {\Bbb N},$
$1\leqslant m<\infty$ such that the following is true: for any
neighborhood $U$ of $x_0$ there is a neighborhood $V\subset U$ of
$x_0$ and  such that:

\medskip
2a) $V\cap D$ is connected,

\medskip
2b) $(V\cap D)\setminus E$ consists at most of $m$ components.

\medskip
Let for $\alpha=2$ the set $F$ have positive capacity, and for
$1<\alpha<2$ it is an arbitrary closed set.

\medskip
Suppose that, for any $x_0\in\partial D$ at least one of the
following conditions $3_1)$--$3_3)$ of Theorem~\ref{th4} is
satisfied.

Let the family of all components of $D^{\,\prime}_f\setminus E_*$ is
equi-uniform over $f\in\frak{R}^{E_*, E, F, N}_{Q, \delta,
\alpha}(D)$ with respect to $\alpha$-modulus. Then every
$f\in\frak{R}^{E_*, E, F, N}_{Q, \delta, \alpha}(D)$ has a
continuous extension to $\partial D$ and the family $\frak{R}^{E_*,
E, F, N}_{Q, \delta, \alpha}(\overline{D}),$ consisting of all
extended mappings $\overline{f}: \overline{D}\rightarrow
\overline{{\Bbb C}},$ is equicontinuous in~$\overline{D}.$ The
equicontinuity must be understood in the sense of the chordal metric
$h$ defined in~(\ref{eq3C}).}
\end{theorem}

\medskip
{\bf Acknowledgements.} The work was supported by the National
Research Foundation of Ukraine (Project ``Analogues of
Carath\'{e}odory and Koebe-Bloch theorems for Orlycz-Sobolev
classes'', Project number 2025.02/0010).

\medskip\medskip
{\bf \noindent Victoria Desyatka} \\
Zhytomyr Ivan Franko State University,  \\
40 Velyka Berdychivs'ka Str., 10 008  Zhytomyr, UKRAINE \\
victoriazehrer@gmail.com\

\medskip
{\bf \noindent Alina Halyts'ka} \\
Zhytomyr Ivan Franko State University,  \\
40 Velyka Berdychivs'ka Str., 10 008  Zhytomyr, UKRAINE \\
alinka16091995@gmail.com

\medskip
{\bf \noindent Evgeny Sevost'yanov} \\
{\bf 1.} Zhytomyr Ivan Franko State University,  \\
40 Velyka Berdychivs'ka Str., 10 008  Zhytomyr, UKRAINE \\
{\bf 2.} Institute of Applied Mathematics and Mechanics\\
of NAS of Ukraine, \\
19 Henerala Batyuka Str., 84 116 Slov'yansk,  UKRAINE\\
esevostyanov2009@gmail.com

\end{document}